\documentstyle[12pt]{article}

\newcommand{\ts}{\otimes} 

\newcommand{\pp}{{\bf P}}
\newcommand{\qq}{{\bf Q}}

\newcommand{\oo}{{\cal O}}
\newcommand{\ii}{{\cal I}}

\newcommand{\ra}{\rightarrow}

\newcommand{\mult}{{\rm mult}}

\newtheorem{theorem}{Theorem}
\newtheorem{lemma}[theorem]{Lemma}
\newtheorem{question}[theorem]{Question}
\newtheorem{example}[theorem]{Example}
\newtheorem{conjecture}[theorem]{Conjecture}
\newtheorem{corollary}[theorem]{Corollary}
\newtheorem{proposition}[theorem]{Proposition}
\newtheorem{definition}[theorem]{Definition}
\newcounter{nmb}

\newcommand{\next}{\addtocounter{nmb}{1}}

\newcommand{\nmbthm}{\renewcommand{\theequation}{\thetheorem.\thenmb}}
\newcommand{\nonmbthm}{\renewcommand{\theequation}{\thetheorem}}

\newenvironment{cor}{\setcounter{nmb}{0}\begin{corollary}}{\end{corollary}}

\begin{document}
\setcounter{section}{-1}
\title{Seshadri constants and the geometry of surfaces}
\author{Michael Nakamaye\footnote{Partially 
supported by NSF Grant DMS 0070190}}
\maketitle

\section{Introduction}
\nonmbthm

In this paper we study how locally defined numerical invariants can carry 
global geometric information about algebraic surfaces.  The local invariants
which we will study are Seshadri constants which measure the positivity of
a line bundle near a point.  

\begin{definition} Suppose $X$ is a smooth projective variety, $x \in X$, and
$A$ an ample line bundle on $X$. Then
$$
\epsilon(x,A) = \inf_{C \ni x}\frac{c_1(A) \cap C}{\mult_x(C)};
$$
here the infimum runs over all integral
curves $C \subset X$ passing through $x$. 
\label{d1}
\end{definition}

\noindent  This numerical definition is equivalent to a more intuitive
geometric definition.  In particular, $\epsilon(x,A)$ is the supremum
of all non--negative rational numbers $\alpha$ such that the linear series
$|nA|$ separates $n\alpha$--jets at $x$ for $n$ sufficiently large and
divisible.  Note that if $L$ is a nef line bundle on $X$ then Definition
\ref{d1} still makes sense and $\epsilon(x,L)$ is defined accordingly.  
When $L$ is nef but not ample,
however, $\epsilon(x,L)$ can take the value zero at some or all points of
$X$.

Given the local nature of this definition, it is surprising that Seshadri
constants can carry global information about a variety $X$.  Indeed, at special
points Seshadri constants rarely carry interesting global information: 
for example, 
suppose $\pi: X \ra Y$ is the blow--up of a smooth surface at a point with
exceptional divisor $E$.  Then for any ample line bundle $A$ on $Y$ and any 
sufficiently large
positive integer $n$, $\epsilon(x,\pi^\ast(nA)(-E)) = 1$ at all
 points $x \in E$.
At a very general point, however, this type of behavior cannot occur as the
numerical equivalence class of any cycle through a very general point moves
to cover the entire variety: recall that a point $x \in X$ is called very
general if $x$ belongs to the complement of a countable union of closed, proper
subvarieties of $X$.  Thus at a very general point one can hope
for Seshadri constants to carry some global information.  

In particular, one 
global property which might be captured via Seshadri constants
is whether or not $X$ admits a dominant
morphism to a variety of smaller dimension.  
Indeed suppose $\pi: X \ra Y$ is a surjective map of projective varieties 
with $\dim(Y) \geq 1$ and
let $A$ be an ample line bundle on $Y$.  If $B$ is an ample line bundle on $X$
and $\eta \in X$ a general point then 
$$
\epsilon(\eta, \pi^\ast(nA) + \epsilon B) \ll 1
$$
for $\epsilon \ll 1$
as one can see by intersecting with a curve $C$ contained in the
fibre $\pi^{-1}(\pi(\eta))$.
Note that the line bundles $\pi^\ast(nA) + \epsilon B$ form
an unbounded family as $n$ grows.  We will show that the converse is 
also true in case $X$ is a surface: namely, 
if there exists an unbounded family of ample line bundles
whose Seshadri constant at a very general point is bounded, then $X$ admits
a dominant morphism to a curve.  

In order to state our result formally, we require a few definitions.  
We let
$$
S(X) = \left\{A \in NE(X) \ts \qq \,\,\mbox{\rm ample}\,|\, 
\epsilon(\eta,A) \leq 1\right\},
$$
where $\eta \in X$ is a very general point and $NE(X)$ denotes divisor classes
modulo numerical equivalence.  
For an ample $\qq$--divisor $A$  we let
$$
m(A) = \sup_{D \equiv A}\left\{\mult_\eta(D)\,|
\,\,D \in \,\,{\rm Div}(X) \ts
\qq \,\,\mbox{effective}\right\}.
$$
To see that the invariant $m(A)$ is well--defined, choose a very ample divisor
$B$ and let $d = \dim(X)$.  
For any $\qq$--divisor $D \equiv A$, we can choose general divisors
$B_1,\ldots,B_{d-1}$, containing $\eta$,  so that the intersection
$$
{\rm support}(D) \cap B_1 \cap \ldots \cap B_{d-1}
$$
is proper.  It then follows that 
$$
\mult_\eta(D) \leq c_1(A) \cap c_1(B)^{d-1}
$$
establishing that the supremum exists in the definition of $m(A)$.

\begin{theorem}
Let $X$ be a smooth projective surface and define 
$$
m(X) = \sup_{A \in S(X)} m(A).
$$
Then $m(X) > 2$ if and only if $X$ admits a
surjective morphism $\phi: X \ra C$ to a curve $C$.  
\label{case}
\end{theorem}

\noindent
Note in particular that Theorem \ref{case} implies that $m(X) > 2$ if
and only if $m(X) = \infty$.  Indeed, whenever $X$ fibres over a variety
$Y$ it is easy, as above, to produce a family of line bundles $\{A_i\}$
with bounded Seshadri constant at a very general point but where 
$m(A_i) \ra \infty$ as $i \ra \infty$.  
One interesting corollary of Theorem \ref{case} is the following:
\begin{cor}
Suppose $A$ is an ample line bundle on a surface
$X$ satisfying 
$$
\sqrt{c_1(A)^2} > \sqrt{3} \epsilon(\eta,A).
$$
Then there exists a non--trivial fibration $\pi:
X \ra C$ such that the general fibre $F_\eta$ is Seshadri--exceptional for
$A$.  
\label{cc}
\end{cor}

\noindent
Note that an inequality like that of Corollary \ref{cc} holds for
any ample line bundle on a smooth surface, namely
$$
\epsilon(\eta,A) \leq \sqrt{c_1(A)^2}.
$$ 
With somewhat more delicate analysis one can also establish the following
generalization of a result proven for abelian varieties in \cite{N}:
\begin{theorem}
Suppose $A$ is an ample line bundle on $X$ satisfying $A^2 > 1$ and
$\epsilon(\eta,A) = 1$.  Then there exists a non--trivial fibration $\pi:
X \ra C$ such that the general fibre $F_\eta$ is smooth and
Seshadri--exceptional for $A$.  
\label{t2}
\end{theorem}

%Seshadri constants have been most often used \cite{B,BS,el,ekl,N}
%in order to determine the local
%geometry of a pair $(X,L)$ consisting of a smooth algebraic variety and
%a nef line bundle.  In other words, they are often used to generate non--zero
%sections of some line bundle with certain local properties at a point 
%$x \in X$.  From this point of view, a Seshadri exceptional subvariety
%is essentially an obstruction, an obstruction to the existence of sections
%with specified jets.  
%In this note, we would like to take a slightly different 
%point of view and try to study how Seshadri exceptional subvarieties can 
%help to determine the global geometry of $X$.  

The methods we use to prove Theorems \ref{case} and \ref{t2}
are very close to those of \cite{el}, namely
we exploit the fact that any curve $C \subset X$ passing through a very general
point moves in a non--trivial family.  
Theorem \ref{t2} also used the Kodaira--Spencer construction of \cite{el}.  

Note that by work of Oguiso \cite{O}, Seshadri constants can not increase
under specialization and thus if $\eta$ is a very general point then 
$\epsilon(\eta,A)$ is the maximal value achieved by the Seshadri constants
$\epsilon(x,A)$ as $x$ varies over all points of $X$. 
Hwang and Keum \cite{HK} study this maximal value, which they denote by 
$\mu(A)$, and, following the methodology of \cite{ekl}, prove that when
$\mu(A)$ is too small relative to the volume of $A$ then $X$ admits
a fibration like that produced in Corollary \ref{cc}, Theorem \ref{t2},
and Corollary \ref{c1}.  The work of Hwang and Keum is very 
closely related to the present paper, analysing the higher dimensional
case and providing several interesting examples. 

\medskip
\noindent
{\it Acknowledgments} It is a real
pleasure to thank the referee for numerous helpful
comments and suggestions which greatly improved the quality of exposition in
this paper.  The referee also kindly shared with me the very interesting
paper of Hwang and Keum \cite{HK} which has not, to my knowledge, been
published.

\section{Main Results}

\nmbthm
\setcounter{theorem}{4}
\noindent
{\bf Proof of Theorem \ref{case}}
One direction of Theorem \ref{case} is trivial.  Namely, suppose that there 
exists a surjective
map $\phi: X \ra C$ to some curve $C$.  Let $\eta \in X$ be a very
general point and let $P = \phi(\eta)$.  We denote by
$F_P$ the fibre of $X$ over $P$.  Then $F_P$ is nef and 
$\epsilon(\eta,F_P) = 0$.  Thus, if $A$ is any ample line bundle on $X$
we have, for $\alpha$ sufficiently small
$$
\epsilon(\eta,\alpha A + nF_P) \leq 1, \,\,\, \forall n \gg 0
$$
and this is clearly a family $S$ of divisors for which $m(S)$ is unbounded.

For the other direction of Theorem \ref{case}, we require the notion of a 
Seshadri exceptional subvariety
\begin{definition}
Suppose $A$ is an ample line bundle on a variety $X$ and $x \in X$.  An
irreducible subvariety $V$ containing $x$ is called Seshadri exceptional at
$x$ relative to $A$ if
$$
\epsilon(x,A) = \frac{\deg_A(V)}{\mult_x(V)}.
$$
\end{definition}
A result of Campana and Peternell \cite{CP} asserts that there always
exists a Seshadri exceptional subvariety.  
In particular, on a surface $X$, if $\epsilon(\eta,A) < \sqrt{A^2}$ then
$X$ itself can not be Seshadri exceptional and thus there must exist
a Seshadri exceptional curve. 

\setcounter{theorem}{2}

Suppose $m(X) > 2$.  Then by definition we can find an ample divisor 
$A \in S(X)$  and an $\alpha > 0$ with 
\next
\begin{eqnarray}
m(A) > 2 + \alpha.
\label{jj}
\end{eqnarray}
Since $\epsilon(\eta,A)
\leq 1$ it follows that $X$ is not Seshadri exceptional for $A$ and thus
there is a Seshadri exceptional curve $C$ satisfying
\next
\begin{eqnarray}
\frac{A \cdot C}{\mult_\eta(C)} \leq 1.
\label{e0}
\end{eqnarray}
Let $m = \mult_\eta(C)$.
Applying \cite{el} Corollary 1.2 gives
\next
\begin{eqnarray}
C^2 \geq m(m-1).
\label{e1}
\end{eqnarray}

We begin by showing 
\next
\begin{eqnarray}
\mbox{$C$ is smooth at $\eta$ and $C^2 = 0$.} 
\label{zzz}
\end{eqnarray}
Suppose first that $m \geq 2$.  We claim that
\next
\begin{eqnarray}
\epsilon(\eta,C) \geq m-1.
\label{e2}
\end{eqnarray}
To prove \ref{e2},  note first 
that $C^2/\mult_\eta(C) \geq m-1$ by \ref{e1} and
for any curve $C^\prime \neq C$ through $\eta$ one has $C \cdot C^\prime/
\mult_\eta(C^\prime) \geq m$.  Thus we have established \ref{e2}.  
Note that since $C^2 > 0$ and $C$ is irreducible, it
follows that $C$ is big and nef.  Consequently, there exists an effective
divisor $D$ so that $C - \delta D$ is ample for any $\delta > 0$.  
Since $\eta$ is general, we can assume that $\eta$ is not contained in $D$.
By \ref{e2}, given $\epsilon > 0$ for all $\delta$ sufficiently small
$$
\epsilon(\eta,C - \delta D) > m-1 - \epsilon.
$$ 
So if $p: Y \ra X$ denotes the blow--up of $X$
at $\eta$ with exceptional divisor $E$ then 
$$
p^\ast(C - \delta D)((-(m-1-\epsilon) E))
$$
is ample.  
Thus for $k$ sufficiently large and divisible the linear series 
$\left|k(C - \delta D) \ts \ii_\eta^{k(m-1-\epsilon)}\right|$ has
an isolated base point at $\eta$ and its general member is irreducible: here
$\ii_\eta \subset \oo_X$ is the ideal sheaf of the point $\eta$.   
Hence for any $\epsilon > 0$ there exists a non--trivial family $\{D_t\}$
of 
$\qq$-divisors, irreducible in a neighborhood of $\eta$, 
numerically equivalent to $C$ with
\next
\begin{eqnarray}
\mult_\eta(D_t) \geq m-1 - \epsilon.
\label{e3}
\end{eqnarray}

By \ref{jj}, there is a $\qq$--divisor $F$, numerically equivalent to A,
with $\mult_\eta(F) >  2+\alpha$.  
Choose a divisor $D_t$ in the family above which meets $F$ properly, except
possibly along $D$.
By \ref{e0} 
$$
A \cdot C \leq m.
$$
On the other hand, using \ref{e3} we have
\begin{eqnarray*}
A \cdot C &=& F \cdot D_t \\
          &\geq& \mult_\eta(F)\mult_\eta(D_t) + \delta D \cdot F  \\
          &\geq& (2+\alpha)(m-1-\epsilon) -O(\delta).
\end{eqnarray*}
For $\delta$ and $\epsilon$ sufficiently small, 
this is a contradiction unless $m = 1$.
We conclude that the first part of \ref{zzz} holds, namely that
$C$  is smooth at $\eta$.  

To prove the second part of \ref{zzz}, note first that $C^2 \geq 0$ since
$C$ passes through a very general point of $X$.  Moreover, since $C$ is
irredubible $C$ must be nef so that $\epsilon(\eta,C)$ is well--defined. 
irreducible
We assume that $C^2 \geq 1$ and derive a contradiction.  
Since $C$ is irreducible we have
$$
\frac{C \cdot C^\prime}{\mult_\eta(C^\prime)} 
\geq 1, \,\,\,\forall \,\, C^\prime \neq C.
$$
Since $C$ is smooth at $\eta$ we also have 
$$
\frac{C \cdot C}{\mult_\eta(C)} \geq 1.
$$
Thus $\epsilon(\eta,C) \geq 1$.  In particular, for $\epsilon > 0$ we
see as in \ref{e3} above that there is a non--trivial
family of divisors $\{D_t\}$, numerically equivalent to $C$ and locally
irreducible at $\eta$, with $\mult_\eta (D_t) > 1 -\epsilon$.  Suppose as above
that $F$ is a $\qq$--divisor numerically equivalent to $A$ with 
$\mult_\eta(F) > 2 + \alpha$.  Arguing exactly as above we see that $F$
cannot intersect $D_t$ properly, a contradiction since the family of
divisors $D_t$ has no base points other than $\eta$.  
Thus we have established the second half of \ref{zzz}, namely that 
$C^2 = 0$,  and we now have a curve which is
a candidate for giving a fibration of $X$.  

\medskip
Since the curve $C$  passes through a very general point $\eta \in X$ 
we can apply
\cite{K} Proposition 2.5 to obtain a non--trivial family of curves in $X$, 
parametrized by a variety $S$, which are numerically equivalent to $C$.
More precisely, there is a scheme $U$ with a flat map $\pi: U \ra S$ to a
reduced scheme $S$ of finite type and a map $f: U \ra X$ such that
$f: \pi^{-1}(s) \ra X$ is  birational for all $s \in S$ and
$C = f(\pi^{-1}(t))$ for some $t \in S$.  
We replace $S$ with a
suitable smooth affine curve $T$
and $U$ with $\pi^{-1}(T)$ so that the new family is also non--trivial.  
Consider the graph $\Sigma \subset U \times T$ of $\pi$ and let
$\Sigma^\prime \subset X \times T$ be the image of $\Sigma$ via the morphism
$f \times {\rm id}$.  
Let $\pi_2:\Sigma^\prime \ra T$ denote the projection
to the second factor. 
Shrinking $T$ if necessary, we can assume that 
$\pi_2^{-1}(t)$ is a curve for all $t \in T$. 
Then for $x \in T$ the divisors
$C_x = \pi_2^\ast(x)$ are algebraically equivalent 
and satisfy $C_x^2 = 0$ since 
each curve $C_x$ is numerically equivalent to $C$.

Choose a general point $x \in T$ and consider the map
$$
\phi: T \ra {\rm Pic}^0(X)
$$
given by $\phi(y) = \oo_X(C_y - C_x)$ where $C_x$ and $C_y$ are the 
curves corresponding to $x,y \in T$.  The map $\phi$ naturally induces a map
$\phi^n: T^n \ra {\rm Pic}^0(X)$ where $T^n$ is the $n^{th}$ Cartesian
product of $T$ and for $n$ sufficiently large $\phi^n$ is not injective.  
This means that for some $m \leq n$ there are 
points $\{P_i\}_{i=1}^m$ and 
$\{Q_i\}_{i=1}^m$, mutually distinct, satisfying
$$
\sum_{i=1}^m (C_{P_i} -C_{Q_i}) \,\,\,\mbox{linearly equivalent to 0}.
$$
Consider the linear system
$$
\left|\sum_{i=1}^m C_{P_i} \right|.
$$
This system contains the two effective divisors, $\sum_{i=1}^m C_{P_i}$ and
$\sum_{i=1}^m C_{Q_i}$.  By hypothesis, these two divisors are disjoint and
thus the linear series $\left|\sum_{i=1}^m C_{P_i} \right|$ is base
point free.  
In particular, using the pencil of divisors spanned by $\sum_{i=1}^m C_{P_i}$
and $\sum_{i=1}^m C_{Q_i}$ gives
 a map $\psi: X \ra \pp^1$ which contracts $C_x$ for all 
$x \in T$ since $\left(\sum_{i=1}^mC_{P_i}\right) \cdot C_x = 0$.  On the
other hand, the map is surjective because the corresponding linear series
has projective dimension one.
 This completes
the proof of Theorem \ref{case}.

\medskip
Note that we have the following immediate corollary to Theorem \ref{case}:

\setcounter{theorem}{5}
\next
\begin{cor}
Suppose $X$ is a smooth projective surface admitting no non--trivial
fibration over a curve.  Then 
$$
\sqrt{A^2} \leq 2 \cdot \epsilon(\eta,A)
$$
for all ample line bundles $A$.  
\label{c1}
\end{cor}
{\bf Proof of Corollary \ref{c1}} Suppose that 
$\sqrt{A^2} > 2 \cdot \epsilon(\eta,A)$.  Then it follows that
there is a ${\bf Q}$--divisor $D$ numerically equivalent to $A$ such that
$\mult_\eta(D) > 2\epsilon(\eta,A)$.  
Arguing as in the proof of Theorem \ref{case} shows
that if $C$ is Seshadri exceptional for $A$ at $\eta$ then $C^2 = 0$ and one 
obtains a fibration of $X$ over a curve.
Corollary \ref{c1} can also be restated as follows: if $A$ is ample on $X$
and $\sqrt{A^2} > 2\epsilon(\eta,A)$ then $X$ fibres over a curve with
general fibre Seshadri exceptional for $A$.  This is slightly weaker than
Corollary \ref{cc} which we prove now.

\medskip

\noindent
{\bf Proof of Corollary \ref{cc}}
According to the proof of Theorem \ref{case}, it suffices
to produce a {\bf Q}--divisor $D$ numerically equivalent to $A$ 
with ${\rm mult}_\eta(D) > 2\epsilon(\eta,A)$.
However, using \cite{ekl} Proposition 2.3 we see that for any $n > 0$
a divisor $D \in |nA|$ with $\mult_\eta(D) > n(\epsilon(\eta,A) + 
\alpha)$ has multiplicity at least $n\alpha$ along $C_\eta$.
Note that it is  critical for this result
 that the  point $\eta$ be very general.
In particular, the cost of imposing order of vanishing $n(\epsilon(\eta,A) +
\alpha)$ at
$\eta$ is asymptotically at most 
$$
\left(\frac{\epsilon(\eta,A)^2}{2} + \alpha \epsilon(\eta,A)\right) n^2,
$$
the maximum cost coming of course when $C_\eta$ is smooth. 
Here $\frac{\epsilon(\eta,A)^2}{2}n^2$ represents the asymptotic cost of
imposing multiplicity $n\epsilon(\eta,A)$ at $\eta$ and $\alpha \epsilon(
\eta,A)n^2$ is the maximal cost of raising the multiplicity by
an additional $n\alpha$.  
Thus one can always obtain a ${\bf Q}$--divisor numerically equivalent
to $A$ with multiplicity at $\eta$
arbitrarily close to 
$$
\epsilon(\eta,A) + \frac{A^2 - \epsilon(\eta,A)^2}{2\epsilon(\eta,A)}. 
$$
In particular, the hypothesis of Corollary \ref{cc} is satisfied
as soon as
$$
\sqrt{A^2} >\sqrt{3} \epsilon(\eta,A).
$$

\setcounter{theorem}{4}

\medskip
\noindent
{\bf Proof of Theorem \ref{t2}}
We begin with  some concrete examples of Corollary \ref{cc}.
Suppose $A^2 \geq 4$ and $\epsilon(\eta,A) = 1$.
Then by Corollary \ref{cc} the Seshadri exceptional curve of $A$ at
$\eta$ must give a 
fibration of $X$ over a curve.  If $A^2 = 3$ and $\epsilon(\eta,A) = 1$ 
the argument for the case where $A^2 \geq 4$
fails if $C_\eta$ is smooth.
In order
to eliminate the possibility that $C_\eta$ is smooth we consider
$C_\eta^2$.  If $C_\eta^2 = 0$ then we obtain the desired fibration.  
If $C_\eta^2 > 0$ then $\epsilon(\eta,C_\eta) \geq 1$ and we find a 
contradiction arguing as in the second half of \ref{zzz}.  

We now consider
the special case where $A^2 = 2$.  If $C_\eta$ is a Seshadri exceptional
curve at $\eta$, then one readily establishes that there
are only three possibilities.
First $C_\eta$ is smooth with $C_\eta^2 = 0$ in which case $C_\eta$ gives
a fibration of $X$ over a curve.  Second, one could have $C_\eta^2 \geq 1$ and
$C_\eta$ smooth at $\eta$.  This would imply that $A \cdot C_\eta = 1$
which contradicts the Hodge index theorem:  
$$
A \cdot C \geq \sqrt{A^2}\cdot \sqrt{C^2}
$$
Finally, one could have $C_\eta^2 = 2$ and $\mult_\eta(C) = 2$:
all other possibilities are eliminated as above using the Hodge index formula 
and \ref{e1}.  
 This last case is
of course possible provided one allows $C_\eta$ to be reducible, namely
one can take $X = C_1 \times C_2$, the product of two curves and 
$A = F_1 + F_2$, the sum of the two fibres through $\eta$.  
In this case, both $F_1$ and $F_2$ 
are Seshadri exceptional and the divisor $A$ satisfies $A^2 = 2$ and 
$\mult_\eta(A) = 2$.  We would now like to establish
that this is essentially the only such possibility.    In particular,
we will show that $C_\eta$ is reducible and thus both components of $C_\eta$
are smooth and give fibrations.

We will assume that $C_\eta$ is irreducible for very general $\eta \in X$ and
derive a contradiction.
We first claim that for a very general point $\eta$ 
we have the numerical equivalence
\next
\begin{eqnarray}
A \equiv C_\eta.
\label{v1}
\end{eqnarray}
Note that $A \cdot C_\eta = 2$ since $\mult_\eta(C_\eta) = 2$ and
$C_\eta$ is Seshadri exceptional for
$A$ at $\eta$.
Thus 
$$
A \cdot (A - C_\eta) = 0.
$$
One checks that $(A-C_\eta)^2 = 0$ and thus the Hodge index theorem implies 
that
$A$ and $C$ are numerically equivalent,
establishing \ref{v1}.
%We will show that for any $\epsilon > 0$ the divisor $A - (1-\epsilon)C_\eta$
%is in the effective cone.  
%Thus $A - C_\eta$ is in the closure of the effective
%cone and $A \cdot (A - C_\eta) = 0$.  By Kleiman's criterion for ampleness,
%this establishes that $A \equiv C_\eta$. 
%Arguing as in the proof of Corollary \ref{cc} and using
%the fact that $\mult_\eta(C_\eta) = 2$ shows that $m(A) \geq 2$.  
%Thus given $\epsilon > 0$  we can choose
%a ${\bf Q}$--divisor $D \equiv A$ such that 
%\next
%\begin{eqnarray}
%\mult_\eta(D) > 2 - \epsilon.
%\label{v2}
%\end{eqnarray}
%We write  $ D = aC_\eta + D^\prime$ where $D^\prime$ intersects $C_\eta$
%properly.  Thus
%\next
%\begin{eqnarray}
%\mult_\eta(D^\prime) = \mult_\eta(D) -2a
%\label{v3}
%\end{eqnarray}
%Combining \ref{v2} and \ref{v3} shows that 
%\begin{eqnarray*}
%2 &=& D \cdot C_\eta \\
%&\geq & aA\cdot C_\eta + i(\eta:D^\prime \cdot C_\eta) \\
%&\geq &   2a + 2(2-\epsilon - 2a).
%\end{eqnarray*}
%Hence $a \geq 1 - \epsilon$, establishing our claim that
%$A - (1-\epsilon)C_\eta$ is in the effective cone.  

Next we claim that the  curve $C_\eta$ with $\mult_\eta(C_\eta) = 2$ and
$A \cdot C_\eta = 2$ is unique.  Suppose that there were two such curves
$C_1$ and $C_2$.  From the previous paragraph, 
we have that $C_1 \equiv C_2 \equiv A$.  Thus $C_1 \cdot C_2 = 2$  but
since both $C_1$ and $C_2$ are singular at $\eta$ this is only possible if
$C_1 = C_2$.  Choose a positive number $m$ so that $mA$ is very ample and
consider the corresponding Chow variety $V$ parametrizing curves $C \subset
X$ of degree $2m$ relative to $A$. Note that any curve $C$ with $mA \cdot C =
2m$ which is singular must either be one of the Seshadri exceptional curves
 $C_\eta$ or must have its singular point in a closed proper subvariety of $X$.
%singular curve in
%the family parametrized by $V$ must be one of the $C_\eta$ if it meets
%an open subset $U \in X$.  
%Let $U \subset X$ be an open subset such that any curve
%of $C$ with $\deg_A(C) = 2$ which is singular at a point in
% $U$ is one of the Seshadri
%exceptional curves $C_\eta$.  
Thus we can find a closed subvariety $W \subset V$ 
parametrizing
curves whose general member is one of the curves $C_\eta$. Let $S \subset
W$ be an open subset parametrizing a family of the curves $\{C_\eta\}$.
%$$
%W = \{v \in V|C_v \cap U \,\,\mbox{is non--empty}\}
%$$
%is open in $V$ and
%$$
%W^\prime = \{w \in W| C_w \cap U \,\,\mbox{is smooth}\}
%$$
%is open we obtain a locally closed subset $S \subset V$ parametrizing the
%curves $C_\eta$.  Since for each $\eta \in U$ there is a unique curve 
%$C_\eta$.
We obtain an embedding $S \hookrightarrow X$ since the curve $C_\eta$ singular
at $\eta$ is unique.  
Consider the universal curve 
$$
D \subset X \times S
$$
defined by the property that $D \cap X \times \eta = C_\eta \times \eta$ 
for all $\eta
\in S$.  We let $F$ denote a local equation defining $D$.  
Note that the divisor $D$ is singular along the diagonal $\Delta$
by Bertini's theorem applied to the map $\pi_2: D \ra S$.

Suppose $x_1,x_2$ are local coordinates on $X$ and $t_1,t_2$ the corresponding
coordinates on $S$. We claim that for all general
 $\eta \in S$ there is a differential operator
$\partial/\partial t$, on the $S$ factor, satisfying
\next
\begin{eqnarray}
\frac{\partial}{\partial t}F|C_\eta = 0,
\label{x1}
\end{eqnarray}
where $(\partial/\partial t)F|C_\eta$ locally represents the Kodaira--Spencer 
class along a one parameter family with tangent direction $t$ at $\eta$.  
Let $x \in C_\eta$ be a general point and choose $\partial/\partial t$ so that
$$
\frac{\partial}{\partial t} F(x) = 0:
$$
this is possible since we have a two parameter 
family of partial derivitives in the coordinates on $S$.  
The Kodaira--Spencer construction in \cite{el} Corollary 1.2 would then 
produce a section of $H^0(C,N)$, where $N$ denotes the normal bundle of $C$
in $X$, vanishing at $\eta$ and at $x$.  Since $\mult_\eta(C_\eta) = 2$ this
would imply that 
$C^2 \geq 3$ unless $\frac{\partial}{\partial t} F|C_\eta = 0$, establishing
\ref{x1}.  

If $A$ denotes an arc in $S$ through $\eta$ with tangent direction $t$ then
the family of curves $\pi_2^{-1}(A)$ has a moving singular point and so as
in \cite{el} \S 2 this  implies that the corresponding
Kodaira--Spencer class is non--trivial, a contradiction.  
We conclude that the curves
$C_\eta$ must be reducible.  The irreducible components of $C_\eta$ give the
fibration in question, establishing Theorem \ref{t2}.

\medskip
\medskip

The question naturally arises whether or not
the fibration in Theorem \ref{case} can be detected in a more intrinsic
fashion by studying the effective cone more closely.
In particular,  one can hope to detect
such a fibration numerically as the fibre $F$ is a nef class satisfying
$c_1(F)^2 = 0$ and, moreover, curves in the equivalence class of $F$
should be Seshadri exceptional
for the appropriate choice of an ample line bundle.  We see, however, in the
following example, that numerical criteria alone are not sufficient
to identify a fibration.

\setcounter{theorem}{6}

\begin{example}  Suppose $E$ is an elliptic curve and consider $X = E \times
E$.  Let $F_1,F_2$ be fibres of the first and second projections respectively
and let $\Delta$ be the diagonal.  Consider the divisor
$$
D = \sqrt{2}F_1 + \sqrt{3}F_2 - \frac{\sqrt{6}}{\sqrt{2}+\sqrt{3}}\Delta.
$$
Then one checks that $D^2 = 0$ and $D \cdot A > 0$ for $A = F_1 + F_2$.  Thus
$D$ is a nef divisor on the boundary of the effective cone. 
On the other hand the ray corresponding to $D$ in the effective cone of $X$
can contain no points corresponding to integral divisor classes on $X$:
to see this, note that 
$$
\frac{D \cdot F_1}{D \cdot \Delta} = \frac{3(\sqrt{3}-\sqrt{2})}{\sqrt{3}+
\sqrt{2}}
$$
is irrational.
Thus the divisor
$D$ is not 
associated to a surjective map $\phi: X \ra C$.  
\label{ee1}
\end{example}
A more natural example than Example \ref{ee1} was given by Mumford (see
\cite{H}).  In particular, Mumford constructs a surface $X$ and
a divisor $D$ wth $D \cdot C > 0$ for all irreducible curves $C \subset X$
but where no multiple of $D$ is effective.  In particular $D$ is nef and
satisfies $c_1(D)^2 = 0$ but $D$ is not associated to a fibration of $X$.

We also have the following phenomenon where a sequence of fibrations,
suitably normalized, can actually converge to a distinct fibration:

\begin{example} Suppose again $X = E \times E$ for an elliptic curve $E$ and
for a postive integer $m$ consider the morphism
$$
\phi_m: E \times E \ra E, \,\,\,\, \phi_m(x,y) = mx - y.
$$
Let $P \in E$ be a point and let $F_m = \phi_m^{-1}(P)$.  Then the divisors
$F_m$ are all on the boundary of the ample cone and $\frac{F_m}{|F_m|}$ 
converges in $NE(X)$ to the class $\frac{F}{|F|}$, $F$ a fibre  of the 
projection to the first factor.
\label{ee2}  
\end{example}

Thus a real nef class $\xi$ with $\xi^2 = 0$ does not necessarily
carry specific geometric information about morphisms from $X$ to a curve.
One still has, however,  the following interesting question
\begin{question} Suppose $X$ is a smooth projective surface
and suppose that there exists a non--zero
nef real class $\xi$ with $\xi^2 = 0$.
Does  $X$ necessarily admit a surjective
morphism $\phi:X \ra C$
to a curve $C$?  
\label{sb}
\end{question}

\noindent
Question \ref{sb} naturally leads to the following

\begin{question} Suppose $X$ is a smooth projective surface and $A$ an ample
line bundle on $X$ with $A^2 = 1$.  If $\eta \in X$ and $\pi: Y \ra X$
the blow--up of $X$ at $\eta$ then does $Y$ admit
 a non--trivial fibration over
a curve $C$ whose general fibre is numerically equivalent to a Seshadri 
exceptional curve of $A$?
\label{c2}
\end{question}
In particular, with the hypotheses of Question \ref{c2},
by \cite{el} we know that 
$\epsilon(\eta,A)
\geq 1$ and hence $\epsilon(\eta,A) = 1$.  Thus the line bundle 
$L = \pi^\ast(A(-E))$ is nef with $L^2 = 0$.  On the other hand, $L$ is
clearly not numerically equivalent to zero so Question \ref{sb} would
provide the desired fibration.

\begin{tabbing}
Department of Mathematics and Statistics\\
University of New Mexico\\
Albuquerque, New Mexico 87131\\
{\em Electronic mail:} nakamaye@math.unm.edu
\end{tabbing}


\begin{thebibliography}{KMM}
\bibitem[B]{B} T. Bauer, {\em Seshadri constants on algebraic surfaces},
    Math. Ann., {\bf 313}, 1999, pp. 547--583.
%\bibitem[BS]{BS} T. Bauer and T. Szemberg, {\em Local positivity of 
%   principally polarized abelian threefolds}, J. reine angew. Math.,
%     {\bf 531}, 2001, pp. 191--200.  
\bibitem[CP]{CP}  F. Campana and T. Peternell, {\em Algebraicity of the 
     ample cone of projective varieties}, J. reine angew. Math., {\bf 
     404}, 1990, pp. 160-166.
%\bibitem[DEL]{DEL} J.--P. Demailly, L. Ein, and R. Lazarsfeld, {\em
%     A subadditivity property of multiplier ideals}, preprint.
\bibitem[EL]{el} L. Ein, R. Lazarsfeld, {\em Seshadri constants on smooth
surfaces} Asterisque {\bf 218}, 1993, pp. 177--186.
\bibitem[EKL]{ekl} L. Ein, O. Kuchle, R. Lazarsfeld, {\em Local positivity
        of ample line bundles}, JDG, {\bf 42}, 1995, pp. 193--219.  
%\bibitem[ELN]{eln} L. Ein, R. Lazarsfeld, and M. Nakamaye, 
%     {\em Zero Estimates, Intersection Theory, and a Theorem of
%     Demailly}, in Andreatta and Peternell eds., {\em Higher Dimensional
%     Complex Varieties}, de Gruyter, 1996, pp. 183--208.
%\bibitem[Fu]{Fu} W. Fulton, {\em Intersection Theory}, Springer, 1984.
%\bibitem[Ha]{Ha} R. Hartshorne, {\em Algebraic Geometry}, Graduate Texts 
%     in Mathematics, Vol. 52, Springer, 1977.
\bibitem[H]{H} R. Hartshorne, {\em Ample subvarieties of algebraic varieties},
     Lecture Notes in Mathematics, {\bf 156}, Springer-Verlag, 1970.
\bibitem[HK]{HK} J--M Hwang and J--H Keum, 
     {\em Seshadri exceptional foliations},
     unpublished manuscript, 2001.
\bibitem[Ko]{K}  J. Koll\'ar, {\em Shafarevich Maps and Automorphic Forms},
     Princeton University Press, 1995. 
%\bibitem[KMM]{KMM}
%     Y. Kawamata, K. Matsuda, and K. Matsuki, {\em Introduction to the
%    Minimal Model program}, In: Oda, T. (ed.) Algebraic Geometry. Proc. Symp.,
%     Sendai, 1985, (Adv. Stud. Pure Math.,{\bf 10}, pp. 283-360).
\bibitem[N]{N} M. Nakamaye, {\em Seshadri constants on abelian varieties}, 
American Journal of Math, {\bf 118}, 1996, pp. 621--635.  
\bibitem[0]{O} K. Oguiso, {\em Seshadri constants in a family of surfaces},
  Math. Ann., {\bf 323}, 2002, pp. 625--631.

\end{thebibliography}
\end{document}